\documentclass{amsart}[12pt]
\usepackage{color}
\pdfoutput=1
\usepackage{amsmath, amsthm, amscd, amsfonts,amssymb}
\usepackage{graphicx,float}
\usepackage[cmtip,arrow]{xy}
\usepackage{pb-diagram,pb-xy}

\newcommand{\boxederrormessage}[1]{\par\noindent\fbox{\vbox{\setlength\parindent{0pt}\noindent#1}}\par}

\DeclareMathOperator{\HLstar}{HL^*}
\newcommand\HL{Hamkins-L\"owe}

\newcommand{\restrict}{{\upharpoonright}}
\newcommand{\cut}{|}

\newcommand{\Rforce}{\mathbb{R}}    % The slecond of the forcings used.
\newcommand{\radin}{\Rforce{}}

\newcommand{\structKC}[1]{(K[C\setminus{#1}],\in,C\setminus{#1})}

\newcommand{\seq}[1]{\langle\,#1\,\rangle}
\newcommand{\set}[1]{\ensuremath{\{ #1 \}}}
\newcommand{\sing}[1]{\set {#1}}

\newcommand{\Top}{\textit{Top}}

\newcommand\KCmodel[1]{\mathcal{K}[#1]}

\DeclareMathOperator{\ps}{\mathcal{P}}

\DeclareMathOperator{\Force}{Force}
\DeclareMathOperator{\Sent}{Sent}

\DeclareMathOperator{\len}{lh}
\DeclareMathOperator{\otp}{otp}
\DeclareMathOperator{\crit}{crit}

\DeclareMathOperator{\col}{Col}

\DeclareMathOperator{\C}{C}

\DeclareMathOperator*{\dintersect}{\triangle}

\DeclareMathOperator{\cof}{cf}

\def\MPB{{\mathbb{P}}}

\DeclareMathOperator{\ult}{Ult} % ultrapower

\newcommand\ufseq{\mathcal U}
\newcommand\extseq{\mathcal E}
\newtheorem{theorem}{Theorem}[section]
\newtheorem{lemma}[theorem]{Lemma}
\newtheorem{proposition}[theorem]{Proposition}
\newtheorem{corollary}[theorem]{Corollary}

\theoremstyle{definition}
\newtheorem{definition}[theorem]{Definition}

\theoremstyle{remark}

\newtheorem{claim}[theorem]{Claim}

\usepackage[pdftex]{hyperref}

\begin{title}[On a question of Hamkins and L\"{owe}]{On a question of Hamkins and L\"{owe} on the modal logic of collapse forcing}
\end{title}

\author[M. Golshani]{Mohammad  Golshani}
\thanks{%
 The first author's research has been supported by a grant from
  IPM (No. 99030417).}

\address{
School of Mathematics, Institute for Research in Fundamental Sciences (IPM), P.O. Box:
19395-5746, Tehran, Iran. }
\email{golshani.m@gmail.com}

\author[W. Mitchell]{ William Mitchell}
\address{Department of Mathematics, University of Florida, 358b
  Little Hall, PO Box 118105, Gainesville, FL 326118105, USA.}
\email{wjm@ufl.edu}

\date{\today}



\begin{document}

\begin{abstract}
Hamkins and L\"{o}we asked whether there can be a model
$N$ of set theory with the property that
$N\equiv N[g]$ whenever $g$
is a generic collapse of a cardinal of $N$ onto~$\omega$.
We give
equiconsistency  results for two weaker versions of this property.
We also include a proof of Woodin's result that the consistency of the
full \HL\ property follows from that of a Woodin cardinal with an
 inaccessible above.
\end{abstract}
 \maketitle

%%%%%%%%%%%%%%%%%%%%%%%%%%%%%%%%%%%%%%%%%%%%%%%%%
%%%%%%%%%%%%%

\section{introduction}
In~\cite{hamkins-lowe}, Hamkins and L\"{o}we introduced the modal
logic of forcing. Suppose $\Gamma$ is a definable class of forcing
notions. A $\Gamma$-translation from propositional modal logic to  set
theory is a map $H: \Sent(\mathcal{L}_m) \rightarrow
\Sent(\mathcal{L}_\in)$ from the set of sentences of the modal logic
to the set of sentences of set theory
which assigns to each propositional variable some sentence in the language of set theory, leaves the propositional connectives unchanged
and in it
 the modal operations $\Diamond$
and $\square$ are interpreted by
\[
H(\Diamond \phi)=\Diamond_\Gamma \phi = \exists \MPB \in \Gamma, \exists p \in \MPB, p \Vdash \phi,
\]
and
\[
H(\square \phi) = \square_\Gamma \phi = \forall \MPB \in \Gamma, \forall p \in \MPB, p \Vdash \phi.
\]
If $\mathcal{H}(\Gamma)$ denotes the set of $\Gamma$-translations, then  we set
\[
\Force(V, \Gamma) = \{\phi \in  \Sent(\mathcal{L}_m) \mid \forall H \in \mathcal{H}(\Gamma), V \models H(\phi) \},
\]
where  $V$ is a model of $ZFC$. We call $\Force(V, \Gamma)$ the modal logic of the class $\Gamma$ over $V$. Also let
\[
\Force(\Gamma) = \bigcap \{ \Force(V, \Gamma)\mid V \models ZFC       \}.
\]
$\Force(\Gamma)$ is called the modal logic of the class $\Gamma$.
The main result of~\cite{hamkins-lowe} says that the modal logic of the class of all forcing notions (over $L$, G\"{o}del's constructible
universe) is exactly $S4.2$, i.e., using the above notation, $\Force(\Gamma) =\Force (L, \Gamma)=S4.2$,
where $\Gamma$ is the class of all forcing notions.

In~\cite{hamkins-leibman-lowe}, among many other things, it is proved that the modal logic of collapse forcing is $S4.3$, i.e.,
$\Force(\col) =S4.3$, where $\col$ denotes the class of all collapsing forcing notions.

Hamkins and  L\"{o}we have asked, in connection with results in~\cite{hamkins-lowe}, whether there can be a model $N$ of $ZFC$  such that $N\equiv N[H]$ whenever $H$ is the generic collapse of any
cardinal onto $\omega$. This question appeared later in Hamkins's paper~\cite{hamkins} as question 10.

Let us explain the impact of the existence of such a model $N$ for the class $\col$ of collapse forcing notions. Recall that a sentence $\phi$
of the language of set theory is called a $\Gamma$-button over $V$ if $V \models \square_\Gamma \Diamond_\Gamma \square_\Gamma \phi$
and it is called a $\Gamma$-switch
over $V$ if $V \models \square_\Gamma \Diamond_\Gamma \phi \wedge \square_\Gamma \Diamond_\Gamma \neg \phi$.
By~\cite{hamkins-leibman-lowe}, if $V$ contains arbitrarily large finite independent families of buttons and switches,
then $\Force(V, \Gamma) \subseteq S5$.
So a model of set theory as $N$ above would be an extreme counterexample in having no
switches at all for the class of collapse forcing, and would have valid principles
of collapse forcing that are beyond $S5$, a hard upper bound for the other natural
classes of forcing.

Modal logics not contained in $S5$ have recently been studied by
Inamdar and L\"{o}we~\cite{inamdar-lowe}, where they determined that
the modal logic of inner models is S4.2\Top, a modal logic which is not contained in
S5. \Top\ is the  axiom

\[
\Diamond(~(\square \phi \leftrightarrow \phi ) \wedge (\square \neg\phi \leftrightarrow \neg\phi )~).
\]
It is clear that if $N$ is a model as above, then we have
\[
N \models \Diamond_{\col}\phi \leftrightarrow \square_{\col}\phi \leftrightarrow \phi
\]
and so clearly
 $\text{\Top} \in \Force(N, \col)$.

\bigskip

The following unpublished theorem of Woodin gives an upper bound for the
consistency of the \HL{} property. With Woodin's permission,
we include a proof in Section~\ref{proof of woodin theorem}.
\begin{theorem}[Woodin]%
\label{thm:main-full}
  If there is a Woodin cardinal with a inaccessible
  cardinal above, then there is a countable well founded model
  satisfying the Hamkins-L\"owe property.
\end{theorem}

We define two weaker variations of the \HL{} property and determine
their consistency strength.  For convenience, we consider the properties for
cardinals up to an inaccessible cardinal $\kappa$ rather than for all
infinite cardinals:
%[[WJM 2/24] _weak_ HL property.]
\begin{definition}%
  \label{def:weakHL}
  Suppose $\kappa$ is an inaccessible cardinal.
  \begin{itemize}\item
    The \emph{Hamkins-L\"owe property on successor cardinals below $\kappa$} holds if $V[g]\equiv V$
    whenever $g\subset\col(\omega,\lambda^+)$ is generic for some
    infinite cardinal $\lambda<\kappa$.
  \item\label{item:HLonClub}
    Suppose that there is no inner model with a Woodin cardinal, and hence the Steel~\cite{mitchell-steel.inner-model} core model $K$ exists. The
    \emph{weak Hamkins-L\"owe property on a club set below $\kappa$} holds if there is a club
    subset $C\subset\kappa$ such that
    \[ \forall\lambda\in C\quad\structKC{\lambda}\equiv
    (K[C], \in, C).\]
  \end{itemize}
\end{definition}

These properties were defined in the course of attempting to determine the exact strength of the \HL\ property.   In Theorem~\ref{thm:main}
we present equiconsistency results---much weaker that a Woodin cardinal---for both weak properties.
We show that in the presence of the Steel core model the \HL\ property
on a club set below $\kappa$ is implied by the \HL\ property below
$\kappa$, and so Theorem~\ref{thm:main} gives the best known
lower bound for the consistency strength of the \HL\ property.

By Tarski's theorem on the undefinability of truth,
using Definition~\ref{def:weakHL}
requires a base theory stronger than ZFC\@.
We will assume as a background theory throughout this paper that
ZFC holds with
an added predicate coding satisfiability for all formulas (without the
added predicate) in the
language.
It should be noted that formulas or sentences refered to in the text never include this extra predicate.

 Before we state our main theorem, let us make the following definition.
\begin{definition}
\label{U-sequence}
If $U$ is a measure on a cardinal $\kappa$ and $E$ is an extender on
$\kappa$ such that
\begin{equation}
  U = E\cut \kappa+1 =\set{x\subseteq\kappa\mid \kappa\in j^{E}(x)},
\end{equation}
then we will call $E$ a $U$-extender.
\end{definition}

%\medbreak\hrule\medbreak
%[[NEW (notes WJM 2/2024)] This replaces the remark after the lemma, which I have deleted.]

The $U$-extender is the first step in the hierarchy of large cardinals beyond sequences of measurable cardinals. Such an extender yields an embedding
\begin{equation}\label{eq:1A}
    j^E\colon V \xrightarrow{j^U} M_1=\ult(V,U)\xrightarrow{k}M
\end{equation}
with $\crit(k)=(\kappa^{++})^{M_1}$. This embedding can be presented as an ultrafilter using the pair $(\kappa,\crit(k))$: If $U'=\set{X\subset \kappa^2\mid (\kappa,\crit(k))\in j(X)}$, then $j^E=j^{U'}$.
% {\color{red} [[Or at least $j^{U'}$ is a $U$-extender subembedding of $j$. FIXME?]]}

The existence of a $U$-extender is stronger than $o(\kappa)=\kappa^{++}$: If $V$ in equation~\ref{eq:1A} is an extender model $L[\extseq]$, then $(\kappa^{++})^{M_1}=\crit(k)\le o^{M_1}(\kappa)$, since $E\in M\setminus \ult(V,E)$.  However it is weaker than $o(\kappa)=\kappa^{++}+1$, that is, an embedding $j\colon V\to M$ with $(\kappa^{++})^M=(\kappa^{++})^V$.

Indeed suppose $V$ is an extender model $L[\extseq]$ and there is an embedding chain of length 3:
\begin{equation}\label{eq:1B}
    j\colon V\xrightarrow{ j^U} M_1=\ult(V,U)\xrightarrow{k_1}M_2\xrightarrow{k_2}M_3
\end{equation}
with $\crit(k_1)=(\kappa^{++})^{M_1}$ and $\crit(k_2)=(\kappa^{++})^{M_2}$, carrying the $U$-extender idea just one step further.
Then there are $\kappa^{++}$ many $U$-extenders, as otherwise all of them would be in $M_2$, but there is at least one which is in $M_3\setminus M_2$: namely that given by the embedding $k_1\circ j^U\colon V\to M_2$.

If $o(\kappa)=\kappa^{++}+1$, then there would be a chain like that of equation~\eqref{eq:1B} of length $\kappa^{++}$: Let $j\colon V\to M$ with $(\kappa^{++})^M=\kappa^{++}$.   Then for every $\alpha<\kappa^{++}$ the ultrafilter $\ufseq_\alpha=\set{x\subset\kappa^2\mid(\kappa,\alpha)\in j(x)}$ factors $j$ into
\begin{equation*}
    j\colon V\xrightarrow{j^{\ufseq_\alpha}}\ult(V,\ufseq_\alpha)\xrightarrow{k}M
\end{equation*}
with $\crit(k)=(\kappa^{++})^{\ult(V,\ufseq_\alpha)}$.  Since $j^{U_\alpha}(\kappa)<\kappa^{++}$, there must be a cofinal set of $\alpha$ which give distinct ultrafilters $\ufseq_\alpha$.

%\medbreak\hrule\medbreak

We will be using
this notion only in extender models $L[\extseq]$; furthermore all extenders we
consider have at most two generators and hence are equivalent to an
ultrafilter.  We will identify $E$ with this ultrafilter.

\begin{theorem}\label{thm:main}
  %Suppose that $\kappa$ is an inaccessible cardinal.
  \begin{enumerate}
  \item  The Hamkins-L\"owe property for successor cardinals below an inaccessible cardinal is
    equiconsistent with the existence of a cardinal
    $\kappa$ such that $o(\kappa)=\kappa^{+}$.
  \item The weak Hamkins-L\"owe property on a club set below an
    inaccessible cardinal is equiconsistent with
    the existence of a cardinal $\kappa$ with a measure $U$ such
    that there are $\kappa^+$  many $U$-extenders.
  \end{enumerate}
\end{theorem}

The following section contains the proof of the right-to-left
direction of Theorem~\ref{thm:main}, establishing the upper bounds for
the consistency strength of the two weak \HL{} properties.
Section~\ref{sec:lower-bound} proves the other direction,
establishing lower bounds for the consistency strength of these two
properties and hence for the full \HL{} property. It also includes the proof that the \HL{} property (assuming no inner model with a Woodin cardinal) implies the weak \HL{} property on a club.

Section~\ref{proof of woodin theorem} gives the proof of Woodin's
Theorem~\ref{thm:main-full}.

\section{Upper bounds for the weak Hamkins-L\"owe properties}
\label{sec:cons-hamk-lowes}
In this section we assume the existence of sequences of ultrafilters
as in Theorem~\ref{thm:main}~(1,2) and force to obtain models with the
weak \HL{} properties.
For the sake of some minor simplifications we assume $V=L[\extseq]$ for some extender sequence $\extseq$.

The primary forcing we use is Radin forcing as introduced in \cite{mitchell82}.
We assume that the reader is acquainted with this forcing and give a
brief definition only for the sake of notation.  A full definition,
and proofs of the properties we use, can be found in~\cite{gitik} or~\cite{mitchell82}.  Note that our definition of coherence for a sequence $\ufseq$ of ultrafilters on $\kappa$ is weaker than that of Radin: it asserts only that for all $\alpha<\len(\ufseq)$ we have $\ufseq\restrict\alpha\in\ult(V,\ufseq_\alpha)$.
%[[FIXME] $\radin_\ufseq$ is the partial order for radin.]

\begin{definition}[ultrafilter sequence]
\label{def:radinforcing}
  An \emph{ultrafilter sequence} on $\kappa$ is a coherent sequence $\ufseq=\ufseq^\kappa$  of $\kappa$ complete ultrafilters on $\kappa$ such that
  \begin{equation*}
    \forall\alpha<\len(\ufseq)\ \ufseq\restrict\alpha\in\ult(V,\ufseq_\alpha).
  \end{equation*}
  For the purpose of this definition, we assume that $\ufseq_0$ is a normal measure on $\kappa$ and the remaining ultrafilters are given in the following isomorphic form:
  \begin{equation*}
    \ufseq_\alpha=\set{x\subset \kappa\mid \ufseq\restrict\alpha\in j^{\ufseq_\alpha}(x)}.
  \end{equation*}
\end{definition}

\begin{definition}[Radin forcing]
\label{def:radinforcing}
 If $\ufseq$ is an ultrafilter sequence on $\kappa$, then the Radin forcing $\radin_\ufseq$
  is defined as follows. A condition in $\radin_\ufseq$ is a pair $(a, A),$ where
  \begin{enumerate}
  \item[(a)] $a \subseteq A \subseteq \kappa,$
   \item[(b)] $a$ is finite,
    \item[(c)] $A \cap \lambda \in U$ whenever $\lambda \in A \cup \{\kappa\}$ and $U$ is a measure on $\lambda.$
  \end{enumerate}
The forcing order on $\radin_\ufseq$ is given by $(a', A') \leq (a, A)$ if and only if $a' \supseteq a, A' \subseteq A$
and $a' \setminus a \subseteq A.$ The order $\leq^*$ is defined as $(a', A') \leq^* (a, A)$ if and only if $a' =a$ and $A' \subseteq A$.
\end{definition}
Note that the forcing $\radin_\ufseq$ is equivalent to the forcing defined in \cite{mitchell82}.
Forcing with  $\radin_\ufseq$ yields a generic object which is a closed unbounded subset of $\kappa$, together with an ultrafilter sequence (possibly empty) $\ufseq^\lambda$ on each $\lambda\in C$ so that for any $\lambda\in C\cup\sing{\kappa}$ and any subset $x\in\ps^V(V_\lambda)$,
  \begin{align}
    \set{\lambda'\in C\cap\lambda\mid \ufseq^{\lambda'}\in x}\quad \text{is cofinal in $\lambda$}\quad &\iff \exists\alpha>0 \ x\in\ufseq_\alpha^\lambda\label{eq:radinCup}\\
    \set{\lambda'\in C\cap\lambda\mid \ufseq^{\lambda'}\notin x}\quad \text{is bounded in $\lambda$}\quad &\iff \forall\alpha>0\ x\in\ufseq_\alpha^\lambda.\label{eq:radinCap}
  \end{align}

Note that the representation of ultrafilters used here is generally not the natural or simplest isomorphic representation of the ultrafilters, but any coherent sequence of ultrafilters can be isomorphically put into this form.   In this paper we will primarily be concerned with the set $C\subset\kappa$.  At the level we are concerned with, the sets $C$ and $G$ are equivalent.

For each clause of Theorem~\ref{thm:main} we start with the Radin
forcing $\radin_{\ufseq}$ where $\ufseq$ is a sequence of ultrafilters of length
$\kappa^+$  as given by the right side of that clause.
Let $\Sigma$ be the set of sentences of the language of set theory
together with a name for the generic set $G$.  Then $\Sigma$ is
countable, and since the direct forcing $\leq^*$ is countably closed,
it follows from the fact that Radin forcing has the Prikry property that there is a condition $p_0=\seq{(\kappa,  B)}$ which
decides every sentence in $\Sigma$.  We will always include the
condition $p_0$ in any generic set.

 No cardinals are collapsed in the generic extension, and
 cofinalities are unchanged except at limit points of $C$.
 In addition, by a theorem in \cite{mitchell82},
 $\kappa$ remains inaccessible since $\ufseq$ has length~$\kappa^{+}$.

\begin{proof}[Proof of Theorem~\ref{thm:main}(2), right to left]
    Let $\ufseq$ be the sequence starting with $\ufseq_0=U$ followed by the first $\kappa^+$ many $U$-extenders in the canonical order in the form specified by Definition~\ref{def:radinforcing}.  Since we are working in $K=L[\extseq]$, this will be a coherent sequence.  Let $C$ be the closed unbounded subset of $\kappa$ corresponding to the generic object $G=\seq{\ufseq^\lambda\mid\lambda\in C}$.

    We want to show that $(K[C],\in,C)\equiv\structKC{\lambda}$ for any $\lambda\in C$, and this will follow if we can show that $C\setminus\lambda$ also comes from a generic $G'\subset \radin_\ufseq$ with $p_0\in G'$.   Taking $G'=\seq{\ufseq^\eta\mid \eta\in C\setminus\lambda}$ would satisfy condition~\eqref{eq:radinCup}, and it also satisfies condition~\ref{eq:radinCup} unless $\lambda$ is a limit member of $C$: in that case, since $(C\setminus\lambda)\cap\lambda$ is empty, that condition requires $\ufseq^{G}$ to be empty.  After $G'$ is altered by setting $\ufseq^{G'}_\lambda$ to the empty sequence, it does satisfy this condition.
\end{proof}

\begin{proof}[Proof of Theorem~\ref{thm:main}(1), right to left]
    Set $\ufseq=\seq{\ufseq_\xi\mid
    \xi<\kappa^{+}} $
    where $\ufseq_{\xi}$ is the  measure on $\kappa$ with
    $o(\ufseq_{\xi})=\xi$, and let $C$ be the
    closed and unbounded subset of $\kappa$ coming from a generic subset
    $G$ of $\radin_{\ufseq}$.
    Then we have
    \begin{equation}
        \label{eq:A}
        \structKC{\lambda}\equiv (K[C],\in, C)
    \end{equation}
    for any successor
    member of $C$.   Let $C^+=C\cup\set{\lambda^+\mid \lambda\in\lim(C)\cap\kappa}$, and let $P^C$ be the Easton support product of the Levy
    collapses $\col(\lambda, {<}\!\min(C\setminus\lambda+1))$ for successor members of $C^+$ as well as for $\lambda=\omega$.
    Let $M=K[C][h]$ where $h\subset P^C$ and $g\subset
    \col(\omega,\min(C))$ are generic.   Thus $C^+$ will be the set of cardinals of $M$ below $\kappa$.

    We claim that $M$
    satisfies the Hamkins-L\"owe property for successor cardinals.
    We need to show that $ M\equiv M^{\col(\omega,\lambda)}$ for any successor member $\lambda$ of $C^+$.
    Let $\lambda$ be a successor member of $C^+$, and let $g\subset\col(\omega,\lambda)$ be the generic collapse.  We want to show that $M\equiv M[g]$, which we will do by writing $M[g]$ as $K[C'][h']$ where $C'=C\setminus\lambda+1$ and $h'$ is a generic subset of $P^{C'}$.   The equivalence will then follow from the homogeneity of $P^{C'}$

    Let $\lambda'=\min(C\setminus\lambda+1)$, the next member of $C$.
    We will let $h'$ be the same as $h$ above $\lambda'$.  Then the missing piece of $h'$ is just the Levy collapse $\col(\omega,{<}\lambda')$, and this needs to absorb the forcing $g$ and the parts of the forcing for $M$ which lie below $\lambda'$.  This forcing is the forcing to obtain $C\cap\lambda'$, followed by $P^{C\cap\lambda'+1}$ and the collapse $g$ of $\lambda$.  Let $\eta=\max(C\cap\lambda)$.  Then $C\cap\eta$ is obtained by Radin forcing at finitely many cardinals less than or equal to $\eta$.  This is a forcing of size at most $\eta^+$, so it has the $\lambda$-chain condition.  The forcings $P^{C\cap\lambda'+1}$ and $\col(\omega,\lambda)$ both have the $\lambda'$-chain condition.
    Thus the possibility of the desired embedding follows from the well known fact that any forcing of size at most $\lambda'$ which requires meeting at most $\lambda'$ dense sets can be embedded into $\col(\omega,\lambda')$.
\end{proof}

Note that this would fail if $\lambda$ were a limit member of $C$.  The proof of the absorption result relies on the fact that the forcing being absorbed has at most $\lambda'$ many dense sets to be met.  If $\lambda$ were a limit member of $C$ then this forcing would involve Ramsey forcing at $\lambda$, which requires meeting $\lambda^+$ many dense sets: one for each set $A\in\bigcap \ufseq^G_\lambda$.  But by the covering lemma, $\lambda'$ must be $(\lambda^+)^K$, since $\lambda$ is singular in $M$.

\section{The strength of a positive answer}
\label{sec:lower-bound}

The aim of this section is to  prove the lower bounds in
Theorem~\ref{thm:main}, that is, to show that the specified weak Hamkins-L\"owe
properties imply the existence in $K$ of the specified ultrafilter
sequences.   However we first give the following result, which both
provides  the
justification for our characterization of \HL\ on a closed unbounded
set as a weak \HL\ property, and, together with
Theorem~\ref{thm:main}, gives a (probably very weak) lower bound for the
strength of the full \HL{} property.

\begin{theorem}[Philip Welch]
  \label{thm:HL2club}
  If there is no inner model with a Woodin cardinal, then the \HL\ property
  implies that there is an unbounded class $C$, closed at all singular
  limit points,  such that
  $(K[C], \in, C)\equiv \structKC{\lambda}$ for all $\lambda\in C$.
\end{theorem}

\begin{proof}
  The key observation is that, since there is no inner model with a Woodin cardinal,
  the core model $K$ exists and for any singular limit
   cardinal $\lambda$ we have
  ${(\lambda^{+})}^{V}={(\lambda^{+})}^{K}$.
  Thus $V^{\col(\omega,\lambda)}$ satisfies that $\omega_1$ is a
  successor cardinal in $K$.  Then the same is true  in
  $V^{\col(\omega,\lambda)}$; that is,  for every  infinite cardinal
  $\lambda$, there is a $K$-cardinal $\eta\in[\lambda,\lambda^{+})$ such
  that $\lambda^{+}={(\eta^{+})}^{K}$.  We will write $\lambda^{*}$ for
  this $\eta$,  and we set $C=\set{\lambda^{*}\mid\omega<\lambda<\kappa}$.

  Then for singular cardinals $\lambda$ we have
  $\lambda=\lambda^{*}$, so $C$  is closed at its singular limit points;
  and whenever $\lambda^{*}\in C$ then
  $\lambda^{*}={(\omega_1^{*})}^{V^{\col(\omega,\lambda)}}$, so
  $C\setminus\lambda^*$ is equal to $C$ as defined in $V^{\col(\omega,\lambda)}$.
\end{proof}
Note that the strength of Theorem~\ref{thm:HL2club} is not from the fact that $C$ is a club set---that is already true for the cardinals of $V$---but that the use of $K\setminus\lambda$ means that a limit cardinal $\lambda$ of $C$ can be the initial member of $C\setminus\lambda$.
\begin{proposition}\label{thm:RadinC}
    Suppose that $\ufseq$ is an ultrafilter sequence on $\kappa$ of length $\kappa^+$.  Then $\ufseq$ and the Radin generic object $\seq{\ufseq^\lambda\mid\lambda\in C}$ can be reconstructed from the club set $C$.
\end{proposition}

\begin{proof}
    Fix $\kappa \leq \alpha< \kappa^+$, and let $r\colon\kappa \leftrightarrow \alpha$ be a bijection in $V$, and for $\lambda\in\kappa$ and $\alpha'\le\alpha$ define $f^r_{\alpha'}(\lambda)=\otp(r[\lambda]\cap \alpha')$.
    Then $[f^r_{\alpha'}]_U=\alpha'$ for all $\alpha'<\alpha$ and all measures $U$ on $\kappa$.

    If $C\subset\kappa$  is $\radin_\ufseq$-generic for an ultrafilter sequence $\ufseq$ on $\kappa$, then we show how to use $f^r_{\alpha}$ to reconstruct $\ufseq$ for $\alpha<\kappa^+$.  For $\lambda$ in a closed unbounded subset of $\lambda$---and hence on a terminal segment of $C$---this process will also use $r'=f^r_{\alpha}\restrict\lambda'$ to find $\ufseq^\lambda\restrict f^r_{\alpha}(\lambda)$.  By applying the same process at $\kappa$ with larger $\alpha$, and by using the same process on initial segments $C\cap\lambda$ of $C$, the full generic sequence can be obtained.

    Define a sequence $\langle D_\nu \mid \nu<\alpha \rangle$ of club subset of $\kappa$ by recursion on $\nu$: $D_0=C$, $D_{\nu+1}=\lim(D_\nu)$, and if $\nu$ is a limit ordinal then
    \begin{equation*}
            D_{\nu}=\dintersect_{\alpha'<\nu}C_{\alpha'}=\set{\lambda<\kappa\mid
            \forall\lambda'<\lambda\; \lambda\in D_{f^r_{\nu}(\lambda')}}.
    \end{equation*}
    Now define $\ufseq_{\kappa,\eta}$, and $\ufseq_{\lambda,f^r_{\eta}(\lambda)}$ for sufficiently large $\lambda\in D_{\eta+1}$, by setting
    $\ufseq_{\kappa,0}=\set{x\subset\kappa\mid D_0\setminus x\text{ is unbounded below }\kappa}$, and
    \begin{equation*}
        \ufseq_{\kappa,\eta}=\set{x\subset V_\kappa\mid D_{\eta}\setminus x\text{ is unbounded below }\kappa},
    \end{equation*}
    and similarly for sufficiently large $\lambda<\kappa$.
\end{proof}

\begin{lemma}
  \label{thm:reversedirection}
  Suppose that $\kappa$ is inaccessible and $C\subset\kappa$ is a  set
  such that $K[C]\equiv K[C\setminus\lambda]$ for all
  $\lambda\in  C$.  Then the method of Lemma~\ref{thm:RadinC} yields an ultrafilter sequence $\ufseq$ on $\kappa$ of length $\kappa^+$ such that $C\cup\lim(C)$ is $\radin_\ufseq$-generic over $K$.
\end{lemma}
\begin{proof}
    Note that if $\sigma\colon C\to\kappa$ is any function defined by a formula
    \begin{equation*}
        \sigma(\nu)=\nu'\iff (K[C \setminus \nu], \in, C \setminus \nu)\models \phi(\nu,\nu')
    \end{equation*}
    where $\phi$ has no other parameters,  if $C$ is closed and unbounded in $\kappa$ and $\sigma(\nu)<\nu$ for some $\nu\in C$, then $\sigma$ is constant on $C$.
   To see this, suppose that $\sigma(\nu)<\nu$ for some $\nu\in C$. Then:
   \begin{itemize}
   \item $\sigma(\nu)<\nu$ for all $\nu\in C$.
   \item $\sigma$ in nondecreasing, i.e.,
    $\sigma(\nu)\leq\sigma(\nu')$ for all $\nu<\nu'$ in~$C$.  If
    not,  then there is $\nu\in C$ such that
    $(K[C\setminus\nu], \in, C\setminus\nu)\models$``$\exists\nu'\in  C \setminus \nu;
    \sigma(\nu')<\sigma(\min(C \setminus \nu))$''.  Thus the same
    sentence is true in $(K[C\setminus\nu], \in, C\setminus\nu)$ for every $\nu\in C$,
    so that for all $\nu\in C$ there is $\nu'>\nu$ in $C$ such that
    $\sigma(\nu)>\sigma(\nu')$,
    but this enables the construction of an infinite decreasing sequence
    of ordinals.
   \end{itemize}

   By Fodor's lemma, $\sigma$ is constant on a stationary set, and since $\sigma$ is nondecreasing, hence it is constant on $C$.
  \iffalse
    \boxederrormessage{03/15/24:
        This also requires $\sigma$ be continuous (if $C$ is continuous.)
    }
    \begin{itemize}
        \item If some $\sigma(\nu)<\nu$, then always $\sigma(\nu)<\nu$.
        \item Nondecreasing.
        \item If ever increasing, then unbounded in $\kappa$.
        \item Then $\sigma[C]$  is like $C$, but with smaller minimal element.
    \end{itemize}
    \fi
    \begin{claim}
        If $\lambda\in\lim(C)$ and $E$ is any closed unbounded subset of $\lambda$, then there is $\lambda'<\lambda$ such that every member of $(\lambda',\lambda)\cap C$ is in $E$.
    \end{claim}
    \begin{proof}[Proof of claim]
        For $\nu\in C$, let $(\lambda,E)$ be the least counterexample in $K[C\setminus\lambda]$, and set $\sigma(\nu)=\max(E\cap\nu)$.   Then $\sigma$ is nondecreasing and   $\sigma(\nu)<\nu$, so $\sigma$ must be constant, contradicting the assumption that $E$ is unbounded in $\lambda$.
    \end{proof}
    Define $D_0=C\cup\lim(C)$, which covers the case when $C$ is not closed, and also includes $\kappa\in D_0$, and define $D_\alpha$ for $\alpha<\kappa^+$ as in the proof of Lemma~\ref{thm:RadinC}.
    Now for $\lambda\in D_0$ and $r\in K$ mapping $\lambda$ onto some $\alpha<\lambda^+$, define $h^r_{\alpha'}$,  $D_{\alpha'}$, and $\ufseq_{\lambda,\alpha'}$ for $\alpha'<\alpha$ as in Theorem~\ref{thm:RadinC}.
    Note that the set of $\nu<\lambda$ such that $r[\nu]\cap\lambda=\nu$ is closed and unbounded in $\lambda$, and hence $h^r_{\alpha'}(\nu)$ is defined for sufficiently large $\nu\in C\cap\lambda$.
    Furthermore, if $r$ and $r'$ are two functions on $\lambda$ mapping $\kappa$ one--one onto $\alpha$ and $\alpha'$ respectively, then $h^r_{\alpha''}(\nu)=h^{r'}_{\alpha''}(\nu)$ for $\alpha''\le\min(\alpha,\alpha')$ and sufficiently large $\nu\in C\cap\lambda$.   Thus the filter $\ufseq_{\lambda,\alpha''}$ is uniquely defined, independently of the choice of $r$.

    Now to show that each $\ufseq_{\lambda,\alpha}$ is a normal, $\lambda$-complete ultrafilter over $K$, suppose that $(\lambda,\alpha,f)$ are minimal counterexamples in $K[C\setminus\nu]$, and set $\sigma(\nu)=f(\nu)$.   Then for $\nu'\in C\cap(\nu,\lambda)$, this is still the least counterexample, so if $\nu'$ is sufficiently large then $\sigma(\nu')<\nu'$.  Thus, $\sigma$ is constant.
\end{proof}

%\medbreak\hrule\medbreak

Clause~(1)  immediately implies the lower bound direction of
Theorem~\ref{thm:main}(1).  Clause~(2) similarly implies the lower bound direction of
Theorem~\ref{thm:main}(2);  in order to slightly simplify the
exposition we make the assumption that $\kappa$ is the smallest
inaccessible limit point of $C$, so that every limit point of $C$ is singular,
and hence $C$ is a closed subset of $\kappa$.
We further assume that $\min(C)$ is as small as possible.
The rest of this section will be devoted to the proof of
Lemma~\ref{thm:reversedirection}.  The proof entails following the construction used in section~(2) in the
proof of lemma~\ref{thm:RadinC}, which showed that the set $C$ obtained from forcing with $\radin_\ufseq$ can be used to reconstruct $\ufseq$ along with the full generic object.   The sets $D_\alpha$ for $\alpha<\kappa^+$ are defined as in that proof (but, for clause~(1), starting with $D_0$ equal to the closure of $C$).
For $\lambda \in C$ let $o^{C}(\lambda)$ be the least ordinal $\alpha$ such that $\lambda$ is not in $D_\alpha.$

For clause~(1) of Lemma~\ref{thm:reversedirection}, we will use the first alternative in the proof of Lemma~\ref{thm:RadinC}, that in which $C$ comes from forcing with $\radin_\ufseq$ where $\ufseq$ is a sequence of normal measures.
We will show that the filters $F^{\lambda}_{\alpha}$ are normal
ultrafilters on $\lambda$ in $K$.
Since we do not yet know that the filters generated by $D_\alpha\setminus D_{\alpha+1}$ are ultrafilters, we will write $F^\lambda_\alpha$ instead of $\ufseq^\lambda_\alpha$.
Assuming that the sequence $\seq
{F^\kappa_\alpha\mid\alpha<\kappa^+}$ is not eventually constant, then we will show that
$o^{K}(F^{\lambda}_{\alpha})\geq\alpha$
for each $\lambda\leq\kappa$ and $\alpha<o^{C}(\lambda)$.
Thus $o^{K}(\lambda)\ge o^{C}(\lambda)$,  completing 
the proof of the first clause of Theorem~\ref{thm:main}.
If the sequence is eventually constant then the proof for clause~(2) of the theorem applies, and we obtain the stronger conclusion.

For the second clause of Theorem~\ref{thm:main}, in which $C$ is closed,
all of the ultrafilters $F^\kappa_\alpha$ will be the same.
We will then write $U$ for
this constant value, and modify the definition of the
filters $F^{\kappa}_{\alpha}$ to obtain a sequence of $U$-extenders; this will complete the proof of Theorem~\ref{thm:main}.

Throughout the rest of this section, we will be working in a language
$\mathcal{L}$ of set theory with an additional
symbol $\dot C$ denoting the set $C$.

\begin{definition}
  \label{def:unif_def} We say that a function $\sigma$ with domain $C$
  is \emph{uniformly definable} if there is a formula $\phi(\mu,x)$, without
  other parameters, such
  that for all $\nu\in C$ and $\mu\in C\setminus\nu$ we have
  $K[C\setminus\nu]\models \phi(\mu,x)\iff x=\sigma(\mu)$.
\iffalse
  \boxederrormessage{03/18/24:
      This definition is strange.  It is actually equivalent to simply ``$x=\sigma(\nu)\iff \KCmodel{C\setminus\nu}\models\phi(\nu,x)$'': it requires $\KCmodel{C\setminus\nu}\models \phi(\mu,x)$ when $\nu=\mu$, and the expression $\phi(\mu,x)$ could perfectly well be ``$\KCmodel{C\setminus\mu}\models\phi'(\mu,x)$'' for some $\phi'$.

      I think I may have put this in thinking it helped with the continuity gap, but it doesn't do anything.
      }
\fi
  If there is such a formula using ordinal parameters, then we say that
  $\sigma$ is \emph{uniformly ordinal definable}.
\end{definition}

\begin{lemma}
    \label{thm:nondecreasing}
    If $\sigma\colon C\to \kappa$ is uniformly definable, then $\sigma$
    is either constant or unbounded in $\kappa$.  Furthermore, if
    $\sigma(\nu)<\nu$ for any $\nu\in C$ then $\sigma$ is constant on $C$.

    If $\sigma$ is uniformly ordinal definable then the set
    $C$ is the union of finitely many intervals, on each of which  $\sigma$ is
    nondecreasing and on the largest of which $\sigma$ satisfies the conclusion of the previous paragraph.
\end{lemma}
\begin{proof}
    First we show if $\sigma$ is uniformly definable then
    $\sigma(\nu)\leq\sigma(\nu')$ for all $\nu<\nu'$ in~$C$.  If
    not,  then there is $\nu\in C$ such that
    $K[C\setminus\nu]\models\exists\nu'\in \dot C\;
    \sigma(\nu')<\sigma(\min(\dot C))$.  Thus the same
    sentence is true in $K[C\setminus\nu]$ for every $\nu\in C$,
    so that for all $\nu\in C$ there is $\nu'>\nu$ in $C$ such that
    $\sigma(\nu)>\sigma(\nu')$,
    but this enables the construction of an infinite decreasing sequence
    of ordinals.

    Now suppose that the uniformly definable function $\sigma$ is bounded in
    $\kappa$, say $\sigma[C]\subseteq \lambda<\kappa$.   Since $\cof(\kappa)>\lambda$, it
    follows that $\sigma$ is constant on $C\setminus\gamma$ for
    sufficiently large $\gamma\in C$.  Then $K[C\setminus\gamma]$ satisfies
    the sentence ``$\sigma$ is constant on $\dot C$'', so this sentence
    is already
    true in $K[C]$.
\iffalse
    \boxederrormessage{03/15/24:
        This next paragraph seems to fix what I thought was a gap.  I'm a bit nervous: its so obvious, it seems like it could be something I thought of and then saw a problem.
        % This next paragraph requires that $\sigma$ be continuous, (for clause~(2) of the main theorem),
        % as otherwise the set $\set{\sigma(\nu)\mid\nu\in C}$ it not closed and hence doesn't ``also satisfy the hypothesis''.

        % I don't see any reason why $\sigma$ would necessarily be continuous.
    }
    %%NEW%%
    \fi
    Finally suppose that $\sigma$ is uniformly definable and
    $\sigma(\nu)<\nu$ for some $\nu\in C$.  Then $\sigma(\nu)<\nu$ for all
    $\nu\in C$.
    Since $C$ is closed and unbounded in $\kappa$, it follows that $\sigma$ is constant on a stationary subset of $\kappa$. Since it it nondecreasing, it is eventually constant and hence is constant.

    % If $\sigma$ is not constant then, by the preceding paragraph,
    % it is cofinal in $\kappa$.   In this case $C'=\set{\sigma(\nu)\mid\nu\in C}$ also satisfies the hypothesis, but it has a smaller
    % minimal member $\sigma(\min(C))$, contradicting the choice of $C$.

    \smallbreak

    It only remains to consider a uniformly ordinal definable function
    $\sigma$.   Suppose that that $C$ is not the
    union of finitely many intervals on which $\sigma$ is
    nondecreasing,
    Let $\phi(\alpha, \nu, x)$ be a  formula, with ordinal parameter
    $\alpha$,  defining $\sigma$ and
    consider the following formula $\psi(\gamma)$:
    \begin{quote}
        The formula $\phi(\gamma,\cdot,\cdot)$ uniformly defines a
        function $\sigma$ on $\dot C$ such that $\dot C$ is not the union
        of finitely many intervals on which $\sigma$ is nondecreasing.
    \end{quote}
    This can be expressed in the language $\mathcal{L}$, so the sentence
    ``$\exists\gamma\;\Psi(\gamma)$'' is true in $K[C\setminus\nu]$
    for all $\nu\in C$.  For each $\nu\in C$ let $\gamma_{\nu}$ be the
    least $\gamma$ such that
    $K[C\setminus\nu]\models\Psi(\gamma)$.

    Now write $\sigma_{\nu}$ for the function defined in
    $K[C\setminus\nu]$ by $\phi(\gamma_\nu, \cdot, \cdot)$, and
    let $\lambda $ be least such that $C\cap\lambda$ is not a finite
    union of intervals on which $\sigma_{\min(C)}$ is nondecreasing.
    Then $\gamma_{\nu}\leq\gamma_{\min(C)}$ for all $\nu\in
    C\cap\lambda$; but the function $\nu\mapsto\gamma_\nu$ is uniformly
    definable, and hence nondecreasing on
    $C$.  Thus $\gamma_{\nu}=\gamma_{\min(C)}$ for all $\nu\in C\cap\lambda$.

    Finally, define $\sigma'(\nu)=\sigma_{\gamma_\nu}(\nu)$.  Then
    $\sigma'$ is uniformly definable and hence is nondecreasing; but this
    is impossible since for $\nu\in C\cap\lambda$ we have
    $\sigma'(\nu)=\sigma_{\gamma_\nu}(\nu)=\sigma_{\min(C)}(\nu)$, which by definition is not nondecreasing.

    \smallbreak{}

    Similarly, if $\sigma$ does not satisfy the conclusion of the first
    paragraph of the lemma on some final segment $C\setminus\nu$ of $C$, then let
    $\gamma$ be the least parameter such that the function $\sigma_\gamma$ defined by $\phi(\gamma,.,.)$
    fails to do so.  Then $\gamma$ is definable in $K[C]$, so the function $\sigma_\gamma$ is uniformly definable and so must satisfy the conclusion of the first paragraph.
\end{proof}

\begin{lemma}
  \label{thm:FalIsNormal}
  For each $\lambda\in \lim(C)$ and each
  $\alpha<o^{C}(\lambda)$,
  $F^{\lambda}_{\alpha}$ is, in $K$,  a normal ultrafilter on $\lambda$.
\end{lemma}
\begin{proof}
  Suppose not.  Uniformly define a function $\sigma$ as follows: for
  each $\nu\in C$ let the pair $(\lambda_{\nu},\alpha_{\nu})$
  be the least counterexample with $\lambda_{\nu}>\nu$, and let $f_{\nu}$ be the
  least function in the order of construction of $K$ witnessing that
  $F^{\lambda_{\nu}}_{\alpha_{\nu}}$ is not normal. Finally,   set
  \begin{equation*}
    \sigma(\nu)=
    \begin{cases}
      f_{\nu}(\nu)&\text{if $\alpha_\nu=0$}\\
        f_{\nu}(\min(C^{\lambda}_{\alpha_\nu}\setminus\nu+1))&\text{otherwise.}
    \end{cases}
  \end{equation*}

We claim that $\sigma$ is constant on $C$.   By
Lemma~\ref{thm:nondecreasing} it is enough to find a single
$\nu\in C$ such that
$\sigma(\nu)<\nu$.   Let
$\lambda=\lambda_{\min(C)}, \alpha=\alpha_{\min(C)}$, and
$f=f_{\min(C)}$.  Then
$\lambda_{\nu}=\lambda$,
$\alpha_\nu=\alpha$  and
$f_{\nu}=f$ for all $\nu\in C\cap \lambda$.
If $\alpha=0$ then the choice of $f$ implies that
$\sigma(\nu)=f_{\nu}(\nu)=f(\nu)<\nu$ for every sufficiently
large $\nu\in C\cap\lambda_{\nu}$.
If $\alpha>0$, on the other hand,  $f(\xi)<\xi$ for every sufficiently large
member of $D_{\alpha}^{\lambda}$.   Fix some such $\xi$.
Since $\alpha>0$, the set $C$ is unbounded in $\xi$ and hence
$C\cap(f(\xi), \xi)$ is not empty, and $\sigma(\nu)=f(\xi)<\nu$ for all
$\nu$ in this set.

This completes the proof that $\sigma$ is constant on $C$.  Now,
continuing the notation from the last paragraph, for each member
$\xi\in D_\alpha^{\lambda}$ we have $f(\xi) = \sigma(\nu)$ for all
sufficiently large members of $C\cap\xi$.  Since $\sigma$ is constant
on $C$ it follows that $f$ is  constant for sufficiently large members of
$D^{\lambda} _{\alpha}$ and hence on a set in $F^{\lambda}_{\alpha}$,
contradicting the choice of $f$.
\end{proof}

\begin{corollary}
  \label{thm:main1}
  If $\lambda\in\lim(C)$ and $\alpha<o^C(\lambda)$ then $o^{K}(F^{\lambda}_{\alpha})\geq\alpha$.
\end{corollary}
\begin{proof}
  The proof is by induction on $\lambda\in \lim(C)$.
  Suppose that $o^{K}(F^{\lambda'}_{\alpha'})\geq \alpha'$
  whenever $\lambda'<\lambda$ and $\alpha'<o^{C}(\lambda')$.

  In order to prove that the same is true of $\lambda$ we use a basic
  consequence of the condensation property of extender models.
  For $\lambda'<\lambda$, $\alpha<\lambda^{+}$
  and $\gamma<\lambda'$ write
  $f^{\lambda,\lambda'}_{\alpha}(\gamma) =
  \otp(f^{\lambda}_{\alpha}[\gamma])$ and set
   $g^{\lambda}_{\alpha}(\lambda')=f^{\lambda,\lambda'}_{\alpha}(\lambda')$.
  Then there is, in $K$, a closed unbounded set $E\subset\lambda$ of $\lambda'$
  such that
  \begin{enumerate}
  \item
    ${[g^{\lambda}_{\alpha}]}_{U} = \alpha$  for any any normal measure
  $U$ on $\lambda$, and
  \item
    $f^{\lambda'}_{g^{\lambda}_{\alpha}(\lambda')}=f^{\lambda,\lambda'}_{\alpha}\restrict \lambda'$.
  \end{enumerate}
  For such $\lambda'$ it follows that
  $C^{\lambda}_{\alpha}\cap\lambda' =
C^{\lambda'}_{g^{\lambda}_\alpha(\lambda')}$.
If $\alpha<o^{C}_{\lambda}$, then $F^{\lambda}_{\alpha}$ concentrates
on $\lambda'\in C^{\lambda}_{\alpha}\cap E$.
Since $C^{\lambda}_{\alpha}\restrict \lambda' =
C^{\lambda'}_{g^{\lambda}_{\alpha}(\lambda')}$, it follows by the
  induction hypothesis that $o^{K}(\lambda')\geq
  g^{\lambda}_{\alpha}(\lambda')$, and since
  ${[g^{\lambda}_{\alpha}]}_{F_{\alpha}^{\lambda}}$ it follows that
  $o^K(F^{\lambda}_{\alpha})\geq\alpha$.
\end{proof}

 This completes the proof of Lemma~\ref{thm:reversedirection}(1) and
 hence of Theorem~\ref{thm:main}(1).

 In order to prove Lemma~\ref{thm:reversedirection}(2), and hence
complete the proof of Theorem~\ref{thm:main}, we assume that the set
$C$ is itself closed and unbounded in $\kappa$, so that
$C^{\lambda}_{\alpha}\subset C$ for all $\lambda$ and $\alpha$.     As was pointed out previously, it
follows that $F^{\lambda}_\alpha=F^{\lambda}_{\alpha'}$ whenever
  $\alpha,\alpha'<o^{C}(\lambda)$.
  For each
$\lambda\in C^{\kappa}_{1}$ let $\ufseq_{\lambda}$ be the constant value of
$F^{\lambda}_{\alpha}$.    The sequence $\seq{\ufseq_\lambda\mid
  \lambda\in C_{1}^{\kappa}}$ is uniformly definable in $K[C]$, so we can modify
  the definition of the filters by defining
  $F
^{\lambda}=\seq{F^{\lambda}_{\alpha}\mid \alpha<o^{C}(\lambda)}$ for
  each $\lambda\in C_1$ by
  \begin{equation*}
    F^{\lambda}_{\alpha}=
    \begin{cases}
      \ufseq_{\lambda}&\text{if $\alpha=0$,}\\
      \set{x\subseteq K_{\lambda}\mid
        (\exists\gamma<\lambda)(\forall\lambda'\in
        D^{\lambda}_{1+\alpha}\!\setminus\gamma)\; \ufseq_{\lambda'}\in x}&\text{otherwise.}
    \end{cases}
  \end{equation*}
Then the same
arguments as before will show that $F^{\lambda}_{\alpha}$ is a
$\ufseq_{\lambda}$-extender whenever $\alpha>0$ and $\lambda\in
C^{\lambda}_{1+\alpha}$.
This completes the proof of Lemma~\ref{thm:reversedirection} and
ofTheorem~\ref{thm:main}.

\section{Woodin's upper bound for the \HL\ property}
\label{proof of woodin theorem}

In this section we give a proof,  motivated by~\cite{friedman-welch-woodin},
 of Woodin's Theorem~\ref{thm:main-full} which states that the consistency of the
 \HL{} property follows from that of a Woodin cardinal with an
 inaccessible cardinal above it.

 For each real $x$ let $M_x$, if it exists, be the minimal transitive
 model of $\text{ZFC}$ containing $x$. Then $M_x= L_{\mu(x)}[x]$ for some countable ordinal
 $\mu(x)$.
Note that $x$ is not supposed to be available as a predicate in $M_x$.

 If $d$ is a Turing degree, then we write $M_d$ for $M_x$,
 where $x$ is some (and hence any) real in $d$.
Given a real $x$ let us denote the Turing degree containing $x$ by ${[x]}_T$.

Now suppose that $\kappa$ is a Woodin cardinal and $\lambda > \kappa$
is an inaccessible cardinal.
Let $\MPB=\col(\omega, \kappa)$ be the Levy collapse forcing for making $\kappa$ countable,
and let $G$ be $\MPB$-generic over $V$. Then by \cite[Corollary 6.12]{neeman}, $\Delta^1_2$-determinacy
holds in  $V[G]$, hence by \cite[Theorem 6.3]{woodin}, $\Sigma^1_2$-determinacy
holds in  $V[G]$.
As $\lambda$ remains inaccessible,
$M_d$ exists in $V[G]$  for each Turing degree $d$ of $V[G]$.
By Martin's cone theorem (see \cite[Theorem 7D.15]{mosch}),
it follows that the theory
of $(M_d ,\in)$ is constant on a cone of Turing degrees, which means that there exists a  Turing degree $d$ such
that $(M_e ,\in) \equiv (M_d, \in)$ for any Turing degree $e \geq_T d$.

We claim that $N=M_d$ satisfies the Hamkins-L\"{o}we property.
To see this,  suppose that $\eta$ is an infinite cardinal of $N$ and $H$ is $\col{(\omega, \eta)}_N$-generic over $N$.
Then there exists a real $y$
such that letting $e= {[y]}_T$, we have $e \geq_T d$ and $N[H]=M_e$.
By our choice of $d$, $M_d \equiv M_e$ and hence $N \equiv N[H]$, as required.

\section{Questions}   % Comments and questions}
\label{sec:conclusion}
The principal  remaining question is, of course, the actual consistency
strength of the \HL{} property, and, in particular, whether the \HL{}
property is consistent with the existence of a core model satisfying
the weak covering lemma.
\iffalse
[[WJM 2/24]  Some improvements to the lower bound for the \HL{} property could be made by strengthening the weak \HL{} property on a club set to state that
\begin{equation*}
    \forall\lambda\in C\quad (K[C]^{\col(\omega,\min(C)},{\in},\C)\equiv
    (K[C]^{\col(\omega,\lambda},{\in},C\setminus\lambda)
\end{equation*}
which also follows from the \HL{} property.   The ultrafilter $\ufseq_\lambda$ would then be definable in $(K[C]^{\col(\omega,\lambda},{\in},C\setminus\lambda)$ as the unique measure on $\lambda$ which is generated by a club subset of $\lambda$.

But this may not be implied by the \HL{} property, since $C\cap\lambda$ need not be definable in $V^{\col(o,\lambda)}$.
]
\fi
We are not able to prove a better lower bound for the \HL{} property
than that given by Theorem~\ref{thm:main}(2).   The apparent next step
would be to prove that all of the ultrafilters (or, equivalently,
extenders with two generators)
$F_{\lambda,\alpha}$ on $\lambda$ are equal, which would imply that
there are $\kappa^{+}$ many 'extenders with three generators, all of
which are equal when restricted to the first two generators.
To use the method of
the proof of  Theorem~\ref{thm:main}(2) would require showing that ultrafilters
$\ufseq_\lambda$ are definable (without
parameters other than $\kappa$ and $C\setminus\lambda$) in
$K{[C]}^{\col(\omega,\lambda)}$, and we have been unable to prove
this.
Note that the forcings used in
this paper are undefinable in $K$, as they use the direct extension condition
$p_0$ which  decides
every sentence in the collapsed model,

\medspace

An interesting intermediate question concerns a strengthening of our
notion of the \HL{} property on a club set:
Let us say that the property $\HLstar$ holds
if there is a club set $C$ such that for all $\lambda\in C$,
$(V^{\col(\omega,\lambda)},C\setminus\lambda)\equiv (V,C)$.

If a Woodin cardinal is necessary for the \HL{} property then this is
likely inconsistent; however if there is no inner model with a Woodin
cardinal then the \HL{} property implies
$K{[C]}^{\col(\omega,\min(C))}\models \HLstar$ where $C$ is the class of
cardinals below $\kappa$.

The property $\HLstar$  appears to be quite strong:
for example, if $C\in M$ witnesses
this property then for any $\lambda, \lambda' \in C$ and any formula
$\sigma(C,\lambda)$ true in $K[C]$ there is a set
$C'\subset\lambda'$ in $K{[C]}^{\col(\omega,\lambda')}$ such that
$K[(C\setminus\lambda')\cup
  C']\models\sigma((C\setminus\lambda')\cup C', \lambda')$.

However it is not clear to us that this would require more large cardinal
strength than that of Theorem~\ref{thm:main}(2).

\bibliographystyle{amsplain}
% \bibliography{Hamkin-Lowe-question}

 %\end{comment}

\end{document}